\documentclass[final]{amsart}
\usepackage{amsmath,amsthm,amsfonts,amssymb,indentfirst,xspace}
\usepackage[notref,notcite]{showkeys}
\usepackage[ps2pdf]{hyperref}
\usepackage[initials]{amsrefs}

\newcommand{\qg}{\eqref{eQGTheta}--\eqref{eQGU}\xspace}
\newcommand{\mqg}{\eqref{eTheta}--\eqref{eU}\xspace}

\newcommand{\del}{\partial}
\newcommand{\lap}{\triangle}
\newcommand{\inv}{^{-1}}

\newcommand{\grad}{\nabla}

\newcommand{\divergence}{\grad \cdot}

\renewcommand{\epsilon}{\varepsilon}
\renewcommand{\leq}{\leqslant}
\renewcommand{\geq}{\geqslant}

\newcommand{\R}{\mathbb{R}}
\newcommand{\Z}{\mathbb{Z}}

\newif\iftextstyle
\textstyletrue
\everydisplay\expandafter{\the\everydisplay\textstylefalse}

\newcommand{\abs}[1]{{\iftextstyle\lvert#1\rvert\else\left\lvert#1\right\rvert\fi}}
\newcommand{\norm}[1]{{\iftextstyle\lVert#1\rVert\else\left\lVert#1\right\rVert\fi}}

\newcommand{\raisesup}[2]{\smash{#1^{#2}_{\vphantom{h}}}\vphantom{#1^{#2}}}

\newcommand{\holderspace}[1]{\raisesup{C}{#1}}
\newcommand{\holdernorm}[2]{\norm{#1}_{\holderspace{#2}}}

\newcommand{\sobolevspace}[1]{\raisesup{\dot{H}}{#1}}

\newcommand{\lp}[1]{\raisesup{L}{\!#1}}
\newcommand{\linf}{\lp{\infty}}
\newcommand{\lpnorm}[2]{\norm{#1}_{\lp{#2}}}
\newcommand{\linfnorm}[1]{\lpnorm{#1}{\infty}}

\newcommand{\besovspace}[3]{\dot{B}^{#1}_{#2,#3}}
\newcommand{\besovnorm}[4]{\norm{#1}_{\besovspace{#2}{#3}{#4}}}

\numberwithin{equation}{section}
\allowdisplaybreaks

\newtheorem{theorem}{Theorem}[section]
\newtheorem*{theorem*}{Theorem}
\newtheorem{lemma}[theorem]{Lemma}

\newtheorem{corollary}[theorem]{Corollary}

\theoremstyle{definition}

\theoremstyle{remark}

\newtheorem*{remark*}{Remark}

\newtheorem{case}{Case}
\begin{document}
\title[Global regularity for a modified quasi-geostrophic equation]{Global regularity for a modified critical dissipative quasi-geostrophic equation}
\author{Peter Constantin}
\author{Gautam Iyer}
\author{Jiahong Wu}
\begin{abstract}
In this paper, we consider the modified quasi-geostrophic equation
\begin{gather*}
\del_t \theta + \left(u \cdot \grad\right) \theta  + \kappa \Lambda^\alpha \theta  = 0\\
u = \Lambda^{\alpha - 1} R^{\perp}\theta.
\end{gather*}
with $\kappa > 0$, $\alpha \in (0,1]$ and $\theta_0 \in \lp{2}(\R^2)$. We remark that the extra $\Lambda^{\alpha - 1}$ is introduced in order to make the scaling invariance of this system similar to the scaling invariance of the critical quasi-geostrophic equations. In this paper, we use Besov space techniques to prove global existence and regularity of strong solutions to this system.
\end{abstract}

\thanks{P.C. acknowledges partial support from NSF grant DMS-0504213. G.I acknowledges partial support from NSF grant DMS-0707920, and thanks the University of Chicago for its hospitality and support.}
\maketitle

\section{Introduction}
The $2$-dimensional quasi-geostrophic equations are
\begin{gather}
\label{eQGTheta} \del_t \theta + \left(u \cdot \grad\right) \theta  + \kappa \Lambda^\alpha \theta  = 0\\
\label{eQGU} u = R^{\perp}\theta
\end{gather}
where $\alpha > 0$, $\kappa \geq 0$, $\Lambda = (-\lap)^{1/2}$ is the Zygmund operator, and
$$
R^{\perp}\theta = \Lambda\inv (-\del_2 \theta, \del \theta).
$$
The case $\alpha = 1$ (termed as the critical case) arises in the geophysical study of rotating fluids \cites{bPedlosky}.

In this paper we consider the following modification of the $2$ dimensional dissipative quasi-geostrophic equation:
\begin{gather}
\label{eTheta} \del_t \theta + \left(u \cdot \grad\right) \theta  + \kappa \Lambda^\alpha \theta  = 0\\
\label{eU} u = \Lambda^{\alpha - 1} R^{\perp}\theta
\end{gather}
We assume $\kappa > 0$ and $\alpha \in (0,1]$.

Note that when $\alpha = 1$ this is the critical dissipative quasi-geostrophic equation. The case of $\alpha=0$ arises when $\theta$ is the vorticity of a two dimensional damped inviscid incompressible fluid \cite{bChorinMarsden}. When $\kappa > 0$, $\alpha \in (0,1)$, the dissipation term is the same as that of the supercritical quasi-geostrophic equation, however the extra $\Lambda^{\alpha - 1}$ in the definition of $u$ makes the drift term $(u\cdot \grad) \theta$ scale the same way as the dissipation $\Lambda^\alpha \theta$. Precisely, equations \mqg are invariant with respect to the scaling $\theta_\epsilon(x, t) = \theta(\epsilon x, \epsilon^\alpha t)$, similar to the scaling invariance of the critical dissipative quasi-geostrophic equation.

Our goal in this paper is to show the global existences of smooth solutions to \mqg with $\lp2$ initial data. For the dissipative quasi-geostrophic equations \qg, this problem has been extensively studied, partly because several authors have emphasized a deep analogy between the $2$-dimensional critical dissipative quasi-geostrophic equations and the $3$-dimensional Navier-Stokes equations. While global existence of the Navier-Stokes equations remains an outstanding open problem in fluid dynamics \cites{bFeff,bConst1}, the global existence of the $2$-dimensional quasi-geostrophic equations was recently settled by Kiselev, Nazarov and Volberg \cite{bKiselevNV} in the periodic case.

Using different techniques, the global existence of smooth solutions to \qg (with $\alpha = 1$) was proved in general $\R^n$ by Caffarelli-Vasseur \cite{bCafarelli}. In the supercritical case ($0 < \alpha < 1$) global existence of smooth solutions is still open. The works  \cites{bConstWuRegularity,bConstWuHolder} have extended the framework of Caffarelli-Vasseur \cite{bCafarelli} to apply in this situation, however two parts of this proof require additional assumptions: H\"older continuity of weak solutions, and smoothness of H\"older continuous solutions. In this paper, we show that both these difficulties can be resolved for the modified equation \mqg. We describe briefly outline this below.\medskip

Following Caffarelli-Vasseur \cite{bCafarelli}, the first step is to show that Leray-Hopf weak solutions to \mqg are in fact $\linf$. Using a level set energy inequality this was shown in \cite{bCafarelli} for general equations of the form \eqref{eTheta}, provided $\alpha = 1$ and $\divergence u = 0$. In the case $0 < \alpha < 1$, the same result has been shown in \cites{bConstWuHolder} for the equations \qg. The latter result directly applies in our situation, and thus Leray-Hopf weak solutions to \mqg are automatically $\linf$.\smallskip

The next step is to show that an $\linf$ Leray-Hopf weak solution of \mqg is also H\"older continuous, with some small exponent $\delta$. For $\alpha = 1$, this has again been shown by Caffarelli-Vasseur \cite{bCafarelli} using a diminishing oscillation result and the natural scaling invariance of the critical quasi-geostrophic equations. The paper \cite{bConstWuHolder} generalizes the diminishing oscillation result in the supercritical case. However the natural scaling of \qg when $0 < \alpha < 1$ will not preserve the BMO norm of $u$, which is required in order to apply the diminishing oscillation result. To circumvent this difficulty, \cite{bConstWuHolder} assumes that $u$ is apriori $\holderspace{1-\alpha}$, which gives the desired control on the BMO norm of $u$ after the appropriate rescaling.

We remark however that the natural scaling of \mqg preserves the BMO norm of $u$ for any $\alpha > 0$. Thus the method of Caffarelli-Vasseur can be applied to show that Leray-Hopf weak $\linf$ solutions of \mqg are actually $\holderspace{\delta}$ for some small $\delta$. However, one can directly deduce this from the work \cite{bConstWuHolder}. Note that equation \eqref{eU} guarantees $u \in \holderspace{1-\alpha}$ provided $\theta \in \linf$ which we know to be true for Leray-Hopf weak solutions. Thus the result of \cite{bConstWuHolder} directly applies in this situation and hence weak solutions of \mqg are automatically H\"older continuous with some small exponent $\delta > 0$.\smallskip

The final step is to show that a Leray-Hopf weak solution which is $\holderspace{\delta}$ is a smooth solution. The paper \cite{bConstWuRegularity} shows this for the supercritical quasi-geostrophic equations provided $\delta > 1 - \alpha$, and that result applies in the present case. Thus the only case that requires special attention is that when $0 < \delta \leq 1 - \alpha$. This is the main theorem of this paper, and the only theorem for which we present the complete proof. Following the method of \cite{bConstWuRegularity}, we essentially show that if a Leray-Hopf weak solution of \mqg is spatially $\besovspace{\delta_1}{p}{\infty}$ for some $\delta_1 \in (0,1)$, then it is actually $\besovspace{\delta'}{p}{\infty}$, where $\delta' = \delta_1 + \min\{\delta_1, \alpha\}$. Successive application of this result will guarantee our weak solution is in fact a classical solution, which can be shown to be smooth via well known methods.\medskip

In the next section, we establish our notational convention, and prove improved regularity of H\"older continuous solutions to \mqg (the main theorem). We only provide a proof for two spatial dimensions, but we remark that the proof goes through almost verbatim in higher dimensions. Finally for completeness, we conclude the paper by stating the required theorems from \cites{bCafarelli, bConstWuRegularity, bConstWuHolder} and using them to deduce smoothness of weak solutions of \mqg.

\section{Improved H\"older regularity}\label{sImprovedBesov}
We recall that $\theta$ is a Leray-Hopf weak solution of \mqg if
$$
\theta \in \linf( [0,\infty), \lp{2}(\R^2)) \cap \lp{2}([0,\infty), \sobolevspace{\frac{\alpha}{2}}(\R^2))
$$
and $\theta$ solves \mqg in the distribution sense.

In this section we will show that if for some $\delta_1 \in (0,1)$, a Leray-Hopf weak solution of \mqg is spatially H\"older continuous with exponent $\delta \in (0,1)$, then it is actually (spatially) H\"older continuous with a better exponent $\delta' = \delta + \frac{1}{2}\min\{\delta, \alpha\}$.

We begin with a brief description of our notation. Let $\{\phi_j \;|\; j \in \Z \}$ be a standard dyadic decomposition of $\R^2$. Namely, for each $j \in \Z$, $\phi_j$ is a Schwartz function with Fourier support (compactly) contained in the annulus $2^{j-1} < \abs{\xi} < 2^{j+1}$ and $\sum_j \hat \phi_j(\xi) = 1$ for $\xi \neq 0$.

We define $\Delta_j$ by $\Delta_j f = \phi_j * f$, $S_j = \sum_{k<j} \Delta_j f$, and the (homogeneous) Besov norm of $f$ by
$$
\besovnorm{f}{s}{p}{q} = \begin{cases}
\left( \sum_j \left( 2^{js} \lpnorm{\Delta_j f}{p} \right)^q \right)^{\frac{1}{q}} & \text{if }q < \infty\\
\sup_j 2^{js} \lpnorm{ \Delta_j f}{p} & \text{if } q = \infty
\end{cases}
$$
and the homogeneous Besov space $\besovspace{s}{p}{q}$ to be the set of all $f$ such that $\besovnorm{f}{s}{p}{q} < \infty$.

We refer the reader to \cite{bConstWuRegularity} for a concise statement of standard embedding theorems, and inequalities we use subsequently. For a more detailed account, and proofs we refer the reader to Stein~\cite{bibLittleStein}*{Chapter 5}, Stein~\cite{bibStein}*{p264}, Schlag~\cite{bibSchlagNotes}, or the classical papers of Taibleson~\cites{bibTaibleson1,bibTaibleson2,bibTaibleson3}.

Finally, we need a lower bound on the (dissipative) term that arises in the process of obtaining $\lp{p}$ estimates of \mqg (see \cite{bibWu}, or Chen, Miao, Zhang~\cite{bibChenMiaoZhang}).

\begin{lemma}\label{lLambdaPositivity}
Let $\alpha \in (0,2)$, and $2 \leq p < \infty$, $j \in \Z$ and $f$ be a tempered distribution on $\R^n$. Then there exists $c = c(n, \alpha, p)$ such that
$$
\int_{\R^n} \abs{\Delta_j f}^{p-2} \Delta_j f \Lambda^\alpha \Delta_j f \geq \frac{2^{\alpha j}}{c} \lpnorm{\Delta_j f}{p}^p
$$
\end{lemma}

We now state and prove the main result of this section.

\begin{theorem}\label{tBetterRegularity}
Suppose $\theta$ is a Leray-Hopf weak solution of \mqg such that for some $\delta > 0$, we have $\theta \in \linf([t_0, t_1], \holderspace{\delta} )$. Then for any $t_0' > t_0$, $\theta \in \linf( [t_0', t_1], \holderspace{\delta'} )$ where $\delta' = \delta + \frac{1}{2}\min\{ \delta, \alpha \}$. 
\end{theorem}
\begin{proof}
Let $p > 2$, and $\delta_1 = (1 - \frac{2}{p}) \delta$. Then
\begin{align*}
\besovnorm{\theta_t}{\delta_1}p\infty &= \sup_j 2^{\delta_1 j} \lpnorm{\Delta_j \theta_t}p\\
    &\leq \sup_j 2^{\delta_1 j} \linfnorm{\Delta_j \theta_t}^{1 - \frac{2}{p}} \lpnorm{\Delta_j \theta_t}{2}^{\frac{2}{p}}\\
    &\leq \holdernorm{\theta_t}{\delta}^{1 - \frac{2}{p}} \lpnorm{\theta_t}{2}^{\frac{2}{p}}
\end{align*}
Thus $\theta \in \linf( [t_0, t_1], \besovspace{\delta_1}p\infty )$. Note that we use the notation $\theta_t$ to denote the function $\theta( \cdot, t)$, and not the time derivative of $\theta$.

Now applying $\Delta_j$ to \eqref{eTheta} gives
\begin{equation}\label{eDeltajTheta}
\del_t \Delta_j \theta + \kappa \Lambda^\alpha \Delta_j \theta = - \Delta_j (u \cdot \grad \theta)
\end{equation}
We know that
\begin{multline*}
\Delta_j (u \cdot \grad \theta) = \sum_{\abs{j-k} \leq 2} \Delta_j\left( S_{k-1} u \cdot \grad \Delta_k \theta\right) + \sum_{\abs{j-k} \leq 2} \Delta_j \left( \Delta_k u \cdot \grad S_{k-1} \theta\right) + \\
+ \sum_{k \geq j-1} \sum_{\abs{k-l} \leq 1} \Delta_j\left( \Delta_k u \cdot \grad \Delta_l \theta\right)
\end{multline*}
Multiplying \eqref{eDeltajTheta} by $p \abs{\Delta_j \theta}^{p-2} \Delta_j \theta$, integrating over $\R^2$ and using Lemma \ref{lLambdaPositivity} gives
\begin{equation}\label{eLpnormDeltajTheta}
\del_t \lpnorm{\Delta_j \theta}{p}^p + \frac{\kappa 2^{\alpha j}}{c} \lpnorm{\Delta_j \theta}{p}^p \leq I_1 + I_2 + I_3
\end{equation}
where
\begin{align*}
I_1 &= -p \sum_{ \abs{j-k} \leq 2} \int \abs{\Delta_j \theta}^{p-2} \Delta_j \theta \cdot \Delta_j \left( S_{k-1} u \cdot \grad \Delta_k \theta \right)\\
I_2 &= -p \sum_{ \abs{j-k} \leq 2} \int \abs{\Delta_j \theta}^{p-2} \Delta_j \theta \cdot \Delta_j \left( \Delta_k u \cdot \grad S_{k-1} \theta \right)\\
I_3 &= -p \sum_{ k \geq j-1} \int \abs{\Delta_j \theta}^{p-2} \Delta_j \theta \cdot \sum_{ \abs{j - l} \leq 1} \Delta_j \left( \Delta_k u \cdot \grad \Delta_l \theta \right)
\end{align*}

We first bound $I_3$ directly using H\"older's and Bernstein's inequalities.
\begin{align}
\nonumber \abs{I_3} &\leq c p \lpnorm{\Delta_j \theta}{p}^{p-1} \lpnorm{\Delta_j \grad \cdot \left( \sum_{k \geq j-1} \sum_{\abs{l - k } \leq 1} \Delta_l u \Delta_k \theta \right)}{p}\\
\label{eI3}    &\leq c p \lpnorm{\Delta_j \theta}{p}^{p-1} 2^j \sum_{k \geq j-1} \sum_{\abs{l-k} \leq 1} \linfnorm{\Delta_l u} \lpnorm{\Delta_k \theta}{p}
\end{align}

Similarly for $I_2$.
\begin{align}
\nonumber \abs{I_2} &\leq c \lpnorm{\Delta_j \theta}{p}^{p-1} \sum_{\abs{j-k} \leq 2} \lpnorm{ \Delta_k u}{p} \lpnorm{\grad S_{k-1} \theta}{\infty}\\
\label{eI2}    &\leq c p \lpnorm{\Delta_j \theta}{p}^{p-1} \sum_{\abs{j-k} \leq 2} \sum_{m \leq k-1} \lpnorm{ \Delta_k u}{p} 2^m \lpnorm{\Delta_m \theta}{\infty}
\end{align}

For $I_1$, we note
\begin{align*}
\sum_{\abs{j-k} \leq 2} \Delta_j \left( S_{k-1} u \cdot \grad \Delta_k \theta \right) &= \sum_{\abs{j-k} \leq 2} \left[ \Delta_j, S_{k-1} u \cdot \grad \right] \Delta_k \theta + \sum_{\abs{j-k} \leq 2} S_{k-1} u \cdot \grad \Delta_j \Delta_k \theta\\
    &= \sum_{\abs{j-k} \leq 2} \left[ \Delta_j, S_{k-1} u \cdot \grad \right] \Delta_k \theta + \sum_{\abs{j-k} \leq 2} S_j u \cdot \grad \Delta_j \Delta_k \theta\\
    &\qquad + \sum_{\abs{j-k} \leq 2} (S_{k-1} u - S_j u) \cdot \grad \Delta_j \Delta_k \theta
\end{align*}
where we use the notation $[A,B]$ to denote the commutator $AB - BA$. Since we know $\sum_{\abs{j-k}\leq 2} \Delta_j \Delta_k = \Delta_j$, we have
$$
I_1 = I_{11} + I_{12} + I_{13}
$$
where
\begin{gather*}
I_{11} = -p \sum_{\abs{j -k} \leq 2} \int \abs{\Delta_j \theta}^{p-2} \Delta_j \theta \cdot \left[\Delta_j, S_{k-1} u \cdot \grad\right] \Delta_k \theta\\
I_{12} = -p \int \abs{\Delta_j \theta}^{p-2} \Delta_j \theta \cdot \left(S_j u \cdot \grad \Delta_j \theta\right)\\
I_{13} = -p \sum_{\abs{j -k} \leq 2} \int \abs{\Delta_j \theta}^{p-2} \Delta_j \theta \cdot \left(\left(S_{k-1} u - S_j u\right) \cdot \grad \Delta_j \Delta_k \theta\right)
\end{gather*}

Note that $u$ (and hence $S_j u$) is divergence free, thus $I_{12} = 0$. We bound $I_{13}$ directly using H\"older's inequality:
\begin{align}
\nonumber \abs{I_{13}} &\leq c p \lpnorm{\Delta_j \theta}{p}^{p-1} \sum_{\abs{j-k} \leq 2} \lpnorm{S_{k-1} u - S_j u}{p} \lpnorm{\grad \Delta_j \theta}{\infty}\\
\label{eI13}    &\leq c p \lpnorm{\Delta_j \theta}{p}^{p-1} 2^{(1 - \delta_1)j} \holdernorm{\theta}{\delta_1} \sum_{\abs{j-k} \leq 2} \lpnorm{\Delta_k u}{p}
\end{align}

We now split the analysis into two cases.

\setcounter{case}{0}
\begin{case}\label{caseDeltaSmall}
$\delta_1 < \alpha$.
\end{case}
In this case, we will show that for any $t_0' > t_0$, $\theta \in \linf( [t_0', t_1], \besovspace{2\delta_1}{p}\infty)$ for any $t > t_0$. After this the theorem will follow using standard embedding theorems about Besov spaces.

We first bound $I_2$, $I_3$ further. The idea is to obtain a $2^{(\alpha - 2 \delta_1)j}$ times norms which are apriori controlled on the right. As we shall see, this doubles the regularity of $\theta$.

From \eqref{eI3} we have
\begin{align*}
\abs{I_3}    &\leq c p \lpnorm{\Delta_j \theta}{p}^{p-1} 2^j \holdernorm{u}{\delta_1 + 1 - \alpha} \sum_{k \geq j-1} 2^{-(\delta_1 + 1 - \alpha)k} \lpnorm{\Delta_k \theta}{p}\\
    &=    c p \lpnorm{\Delta_j \theta}{p}^{p-1} 2^{(\alpha - 2\delta_1)j} \holdernorm{u}{\delta_1 + 1 - \alpha} \sum_{k \geq j-1} 2^{(1 + 2\delta_1 - \alpha)(j-k)} 2^{\delta_1 k} \lpnorm{\Delta_k \theta}{p}\\
    &\leq c p \lpnorm{\Delta_j \theta}{p}^{p-1} 2^{(\alpha - 2\delta_1)j} \holdernorm{\theta}{\delta_1} \besovnorm{\theta}{\delta_1}{p}{\infty}.
\end{align*}

For $I_2$, we have from \eqref{eI2}
\begin{align*}
\abs{I_2}    &= cp \lpnorm{\Delta_j \theta}{p}^{p-1} \sum_{\abs{j-k} \leq 2} \lpnorm{ \Delta_k u}{p} 2^{(1 - \delta_1)k} \sum_{m \leq k-1} 2^{(m-k)(1-\delta_1)} 2^{m \delta_1} \lpnorm{\Delta_m \theta}{\infty}\\
    &\leq cp \lpnorm{\Delta_j \theta}{p}^{p-1} \holdernorm{\theta}{\delta_1} 2^{(\alpha - 2\delta_1)j} \sum_{\abs{j-k} \leq 2}  2^{(k-j)(\alpha - 2\delta_1)} 2^{(\delta_1 + 1 - \alpha)k} \lpnorm{ \Delta_k u}{p}\\
    &\leq cp \lpnorm{\Delta_j \theta}{p}^{p-1} 2^{(\alpha - 2\delta_1)j} \holdernorm{\theta}{\delta_1} \besovnorm{u}{\delta_1 + 1 - \alpha}{p}{\infty}\\
    &\leq cp \lpnorm{\Delta_j \theta}{p}^{p-1} 2^{(\alpha - 2\delta_1)j} \holdernorm{\theta}{\delta_1} \besovnorm{\theta}{\delta_1}{p}{\infty}
\end{align*}

For $I_1$, we bound $I_{11}, \dots, I_{13}$ individually. For $I_{13}$ we have from \eqref{eI13}
\begin{align*}
    &= c p \lpnorm{\Delta_j \theta}{p}^{p-1} 2^{(\alpha - 2\delta_1)j} \holdernorm{\theta}{\delta_1} \sum_{\abs{j-k} \leq 2} 2^{(j-k)(\delta_1 + 1 - \alpha)} 2^{(\delta_1 + 1 - \alpha)k} \lpnorm{\Delta_k u}{p}\\
    &\leq c p \lpnorm{\Delta_j \theta}{p}^{p-1} 2^{(\alpha - 2\delta_1)j} \holdernorm{\theta}{\delta_1} \besovnorm{\theta}{\delta}{p}{\infty}
\end{align*}

The term $I_{12} = 0$ and requires no bounding. Finally we bound the commutator $I_{11}$. Note that
$$
\left[\Delta_j, S_{k-1} u \cdot \grad\right] \Delta_k \theta = \int \phi_j(x-y) \left[S_{k-1}u(y) - S_{k-1}u(x)\right] \cdot \grad \Delta_k \theta(y) \, dy
$$
Since $\delta_1 < \alpha$, $\delta_1 + 1 - \alpha < 1$, thus
\begin{align*}
\linfnorm{ S_{k-1} u(x) - S_{k-1}u (y)} &\leq \holdernorm{u}{\delta_1 + 1 - \alpha} \abs{x - y}^{\delta_1 + 1 - \alpha}\\
    &\leq c \holdernorm{\theta}{\delta_1} \abs{x - y}^{\delta_1 + 1 - \alpha}.
\end{align*}
Hence
\begin{align*}
\abs{I_{11}} &\leq c p \lpnorm{\Delta_j \theta}{p}^{p-1} 2^{-(\delta_1 + 1 - \alpha)j} \holdernorm{\theta}{\delta_1} \sum_{\abs{j-k} \leq 2} 2^k \lpnorm{\Delta_k \theta}{p}\\
    &\leq c p \lpnorm{\Delta_j \theta}{p}^{p-1} 2^{(\alpha - 2\delta_1)j} \holdernorm{\theta}{\delta_1} \sum_{\abs{j-k} \leq 2} 2^{\delta_1 k} \lpnorm{\Delta_k \theta}{p}\\
    &\leq c p \lpnorm{\Delta_j \theta}{p}^{p-1} 2^{(\alpha - 2\delta_1)j} \holdernorm{\theta}{\delta_1} \besovnorm{\theta}{\delta_1}{p}{\infty}
\end{align*}

Combining estimates, we have from \eqref{eLpnormDeltajTheta}
\begin{equation}\label{eLpnormDeltajTheta1}
\del_t \lpnorm{\Delta_j \theta}{p} + \frac{\kappa 2^{\alpha j}}{c} \lpnorm{\Delta_j \theta}{p} \leq c 2^{(\alpha - 2 \delta_1)j} \holdernorm{\theta}{\delta_1} \besovnorm{\theta}{\delta_1}{p}{\infty}
\end{equation}
which upon integration yields
$$
\lpnorm{ \Delta_j \theta_t }{p} \leq e^{-\frac{\kappa 2^{\alpha j} }{c} \left(t - t_0\right)} \lpnorm{\Delta_j \theta_{t_0}}{p} + c\int_{t_0}^t e^{-\frac{\kappa 2^{\alpha j} }{c} (t - s)} 2^{(\alpha - 2 \delta_1)j} \holdernorm{\theta_s}{\delta_1} \besovnorm{\theta_s}{\delta_1}{p}{\infty} \, ds.
$$
Multiplying by $2^{2\delta_1 j}$ and taking the supremum in $j$ gives
\begin{multline*}
\besovnorm{\theta_t}{2\delta_1}{p}{\infty} \leq \sup_j e^{-\frac{\kappa 2^{\alpha j} }{c} \left(t - t_0\right)}2^{2\delta_1 j} \lpnorm{\Delta_j \theta_{t_0}}{p} +\\
+ \frac{c}{\kappa} \sup_j \left(1 - e^{- \frac{\kappa 2^{\alpha j}}{c} \left(t - t_0\right)} \right) \sup_{s \in [t_0,t]} \holdernorm{\theta_s}{\delta_1} \besovnorm{\theta_s}{\delta_1}{p}{\infty}
\end{multline*}
which immediately shows that for any $t_0' > t_0$, $\theta \in \linf( [t_0', t_1], \besovspace{2\delta_1}p\infty )$.

Now note that
$$
2\delta_1 - \frac{2}{p} = 2\left( \delta - \frac{2}{p} \right) - \frac{2}{p}
$$
and hence as $p \to \infty$, $2\delta_1 - \frac{2}{p} \to 2 \delta$. Thus for some large choice of $p$, we have $2 \delta_1 - \frac{2}{p} = \frac{3\delta}{2}$.
Thus for this $p$, we have
$$
\besovspace{2\delta_1}p\infty \subset \besovspace{3\delta/2}\infty \infty
$$
by the Besov embedding theorem. Finally we know $\linf \cap \besovspace{3\delta/2}\infty \infty = \holderspace{\frac{3\delta}{2}}$, concluding the proof for Case \ref{caseDeltaSmall}.
\begin{case}\label{caseDeltaBig}
$\delta_1 \geq \alpha$.
\end{case}
This case can already be handled by result of \cite{bConstWuRegularity}, and we only provide a brief sketch here for completeness. The main difference here is in the commutator $I_{11}$, where we can only get a $2^{-\delta_1 j}$ on the right. Consequently, this will increase the regularity of $\theta$ by $\alpha$ (and not $\delta_1$, as in the previous case).


We deal with the commutator $I_{11}$ first. Note that $\delta_1 \geq \alpha$ implies $\delta_1 + 1 - \alpha \geq 1$, and hence
\begin{align*}
\linfnorm{ S_{k-1}u(x) - S_{k-1} u(y) } &\leq \linfnorm{\grad u} \abs{x - y}\\
    &\leq \holdernorm{\theta}{\delta_1} \abs{x - y}
\end{align*}
This in turn gives
\begin{align*}
\abs{I_{11}} &\leq c p \lpnorm{\Delta_j \theta}{p}^{p-1} 2^{-j} \holdernorm{\theta}{\delta_1} \sum_{\abs{j-k} \leq 2} 2^k \lpnorm{\Delta_k \theta}{p}\\
    &\leq c p \lpnorm{\Delta_j \theta}{p}^{p-1} 2^{-\delta_1 j} \holdernorm{\theta}{\delta_1} \sum_{\abs{j-k} \leq 2} 2^{\delta_1 k} \lpnorm{\Delta_k \theta}{p}\\
    &\leq c p \lpnorm{\Delta_j \theta}{p}^{p-1} 2^{-\delta_1 j} \holdernorm{\theta}{\delta_1} \besovnorm{\theta}{\delta_1}{p}{\infty}
\end{align*}

The bounds for $I_2$, $I_3$ and $I_{13}$ are similar to the first case, and we omit the details. Combining our estimates leads us to \eqref{eLpnormDeltajTheta1} with $2^{(\alpha - 2\delta_1)j}$ replaced with $2^{-\delta_1 j}$. Multiplying by $2^{(\alpha + \delta_1)j}$ and integrating gives
\begin{multline*}
\besovnorm{\theta_t}{\delta_1 + \alpha}{p}{\infty} \leq \sup_j e^{-\frac{\kappa 2^{\alpha j} }{c} \left(t - t_0\right)}2^{(\alpha + \delta_1) j} \lpnorm{\Delta_j \theta_{t_0}}{p} +\\
+ \frac{c}{\kappa} \sup_j \left(1 - e^{- \frac{\kappa 2^{\alpha j}}{c} \left(t - t_0\right)} \right) \sup_{s \in [t_0,t]} \holdernorm{\theta_s}{\delta_1} \besovnorm{\theta_s}{\delta_1}{p}{\infty}.
\end{multline*}
As before, this shows that for any $t_0' > t_0$, $\theta \in \linf( [t_0', t_1], \besovspace{\delta_1 + \alpha}{p}\infty )$.

Now, $\delta_1 + \alpha - \frac{2}{p}$ converges to $\delta + \alpha$ as $p \to \infty$. Thus for some large $p$, we must have $\delta_1 + \alpha - \frac{2}{p} = \delta + \frac{\alpha}{2}$. Applying the Besov embedding concludes the proof in Case \ref{caseDeltaBig}.
\end{proof}
\section{Regularity of weak solutions}

Given Theorem \ref{tBetterRegularity}, one can use the work \cites{bConstWuHolder} and \cite{bCafarelli} to immediately show the existence of global smooth solutions to \mqg with $\lp2$ initial data. We recall the relevant facts from \cites{bConstWuHolder,bConstWuRegularity,bCafarelli} in this section, and briefly outline the proof.

\begin{theorem}[Caffarelli-Vasseur \cite{bCafarelli}, Constantin-Wu \cite{bConstWuHolder}]\label{tL2Linf}
Let $\theta_0 \in \lp{2}(\R^2)$, and $\theta$ be a Leray-Hopf weak solution of \mqg with initial data $\theta$. Then for any $t > 0$, $\theta_t \in \linf(\R^2)$, and further
$$
\linfnorm{\theta_t} \leq c \frac{\lpnorm{\theta_0}{2}}{(\kappa t)^{1/\alpha}}
$$
\end{theorem}

We remark that Caffarelli-Vasseur \cite{bCafarelli} only proves Theorem \ref{tL2Linf} for $\alpha = 1$, and Constantin-Wu \cite{bConstWuHolder} only prove Theorem \ref{tL2Linf} for the system \qg. The proof of this theorem in Constantin-Wu \cite{bConstWuHolder} however only uses the fact that $u$ is divergence free, and thus applies directly for the system \mqg. We do not present the proof of Theorem \ref{tL2Linf} here.

\begin{corollary}\label{cuHolder}
Under the assumptions of Theorem \ref{tL2Linf}, for any $t > 0$, $u_t \in \holderspace{1-\alpha}$ and further
$$
\holdernorm{u_t}{1-\alpha} \leq c \frac{\lpnorm{\theta_0}{2}}{(\kappa t)^{1/\alpha}}
$$
\end{corollary}
\begin{proof}
This follows immediately from the fact that
\begin{equation*}
\holdernorm{ \Lambda^{\alpha-1} f }{1-\alpha} \leq c \linfnorm{f}\qedhere
\end{equation*}
\end{proof}
\begin{corollary}\label{cHolder}
Under the assumptions of Theorem \ref{tL2Linf}, for any $t_0 > 0$, $\theta \in \holderspace{\delta}(\R^2 \times [t_0,\infty))$ for some $\delta > 0$.
\end{corollary}
\begin{proof}
By Corollary \ref{cuHolder}, we know $u \in \linf( [t_0,\infty), \holderspace{1-\alpha}(\R^2) )$. Thus the results of Constantin and Wu \cite{bConstWuHolder} (Theorem 4.1 in particular) applies proving the corollary.
\end{proof}

\begin{lemma}
Suppose $\theta$ is a Leray-Hopf weak solution of \mqg. If for any $\theta \in L^\infty( [t_0, t_1], C^\delta(\R^2)$ for some $\delta \in (0,1)$, then $\theta \in C^\infty( (t_0, t_1] \times \R^2 )$.
\end{lemma}
\begin{proof}
We apply Theorem \ref{tBetterRegularity} can be applied repeatedly to show that for any $t_0' > t_0$, $\theta \in \linf( [t_0', t_1], \holderspace{\delta'} )$ for some $\delta' > 1$. Now the space regularity can be converted to time regularity, showing that $\theta$ is a classical solution of \mqg on the interval $[t_0', t_1]$. Higher regularity now follows via standard techniques.
\end{proof}

\begin{theorem}\label{tSmooth}
For any $\theta_0 \in \lp2(\R^2)$, there exists $\theta \in C^\infty( \R^2 \times (0,\infty))$ which solves \mqg with initial data $\theta_0$.
\end{theorem}
\begin{proof}
Global existence of Leray-Hopf weak solutions to \mqg can be established using the standard method of Galerkin approximations (see for instance \cite{bResnick}, in the case of \qg, or \cite{bConstFoias} in the case of Navier-Stokes). The proof is now immediate from the above results.
\end{proof}

\section*{Acknowledgement}

Stimulating discussions with Luis Caffarelli are gratefully acknowledged.

\end{document}